\documentclass[a4paper,11pt]{article}
\usepackage[cp1251]{inputenc}
\usepackage[english]{babel}
\usepackage{amsmath}
\usepackage{mathrsfs}
\usepackage{amsfonts}
\usepackage[all]{xy}
\usepackage{amscd}
\usepackage{amssymb}

\begin{document}

\title{Trivial $I^\tau$ fibrations of the multiplication maps for monads $\mathbb{O}$ and
$\mathbb{OH}$}

\author{Lesya Karchevska}
\date{}
\maketitle

\begin{abstract}
In this paper we further investigate the geometry of monads
$\mathbb{O}$ of order-preserving functionals and $\mathbb{OH}$ of
positively homogeneous functionals. We prove that for any $X\in
Comp$ with $w(X) = \tau$ the map $\mu_F X$, where $F\in\{O,OH\}$, is
homeomorphic to trivial $I^\tau$-fibration if and only if $X$ is
openly generated $\chi$-homogeneous compactum.
\end{abstract}

\medskip  Department of Mechanics and Mathematics, Ivan Franko National University of Lviv,
Universytetska st.,1 79602 Lviv Ukraine, e-mail:
\emph{crazymaths@ukr.net}.



\vspace*{\parindent} \vspace*{\parindent}

{\bf Introduction.} The geometric properties of different functors
have been studied extensively over the past few decades
\cite{CovFunkt}. The researches concern studying the question of how
functors affect properties of spaces and maps between them as well
as the investigation of properties of maps involved in the
structures generated by functors (i.e. monad multiplication maps,
structural mappings of algebras).

In the present paper we continue the investigation of geometric
properties of monads which have functional representation, that is,
which can be embedded in monad $\mathbb{V}$. This research concerns
the monads $\mathbb{O}$ and $\mathbb{OH}$ generated by functors of
order-preserving and positively homogeneous functionals respectively
and is to answer the question when multiplication maps for these
monads are trivial fibrations with fibers homeomorphic to the
Tychonov cube.

The results of M. Zarichnyi on the inclusion hyperspaces monad (see
\cite{ZarMonMult}), for instance, show that for a continuum $X$ the
multiplication map $\mu_G X$ for this monad is homeomorphic to the
projection map $pr_{G(X)}:I^\tau\times G(X)\to G(X)$ iff $X$ is
openly generated and $\chi$-homogeneous. In this research we obtain
a similar condition for multiplication maps of monads $\mathbb{O}$
and $\mathbb{OH}$ ($X$ is not necessarily connected in our case). It
is interesting whether there is some kind of a general result of
this type (in case of submonads of $\mathbb{V}$, for example).

\vspace*{\parindent}

{\bf Definitions and facts.} In this section we shall recall some
necessary definitions and results from infinite-dimensional topology
as well as define the objects of our investigation - monads of
order-preserving and positively homogeneous functionals and name
some of their properties.

Since in what follows we will deal with endofunctors in the category
$Comp$, we assume all spaces to be compact Hausdorff (briefly
compacta) and mappings to be continuous.

By $w(X)$ we denote the weight of a space $X$, and by $\chi(x,X)$
the character at a point $x\in X$. We call $X$ $\chi$-homogeneous if
for every $x,y\in X$ we have $\chi(x,X) = \chi(y,X)$.

We say that a space $X$ is a \emph{retract} of $Y$, where $X\subset
Y$, if there exists a map $r:Y\to X$ with $r\vert_X =
\mathrm{id}_X$. The space $X$ is called an \emph{absolute retract}
(briefly an $AR$), if for every embedding $i:X\hookrightarrow Y$ the
subspace $i(X)$ is a retract of $Y$.

Recall that a \emph{$\tau$- system}, where $\tau $ is any cardinal
number, is a continuous inverse system consisting of compacta of
weight $\le \tau$ and epimorphisms over a $\tau$-complete indexing
set. As usual, $\omega$ stands for the countable cardinal number. A
compactum $X$ is called \emph{openly generated}, if it can be
represented as the limit of some $\omega$-system with open bonding
mappings \cite{Shchepin}.

By $C(X)$, where $X\in Comp$, we denote the Banach space of all
continuous real-valued functions on $X$ with the sup-norm
$\Vert\varphi\Vert=\sup\lbrace\vert\varphi(x)\vert \ \vert \ x\in
X\rbrace$. By $c_X$, where $c\in \mathbb{R}$, we denote the constant
function: $c_X(x)=c$ for all $x\in X$.

Let $\nu :C(X)\to \mathbb{R}$ be a functional. We say that $\nu$ is
1) normed, if $\nu(1_X)=1$; 2) weakly additive, if for any $\phi \in
C(X)$ and $c\in \mathbb{R}$ we have $\nu(\phi+c_X)=\nu(\phi)+c$; 3)
order-preserving, whenever for any $\varphi,\psi\in C(X)$ such that
$\varphi(x)\le\psi(x)$ for all $x\in X$ (i.e. $\varphi\le \psi$) the
inequality $\nu(\varphi)\le\nu(\psi)$ holds; 4) positively
homogeneous, if for any $\varphi\in C(X)$ and any real $t\ge 0$ we
have $\nu(t\varphi)=t\nu(\varphi)$.

Now for any $X\in Comp$ denote $V(X) = \prod_{\varphi\in
C(X)}[\min\varphi,\max\varphi]$. For any mapping $f:X\to Y$ let
$V(f)$ be mapping such that $V(f)(\nu)(\varphi) = \nu(\varphi\circ
f)$ for any $\nu\in V(X), \ \varphi\in C(Y)$. Defined that way, $V$
forms a covariant functor in category $Comp$.

For any space $X$ by $O(X)$ denote the set of functionals satisfying
1)--3) (order-preserving functionals), by $OH(X)$ -- the set of all
functionals on $C(X)$, which satisfy properties 1)--4) (positively
homogenous functionals). Let $F$ stand for one of $O,OH$. The space
$F(X)$ is considered as the subspace of $V(X)$. For any function
$f:X\to Y$ the map $F(f):F(X)\to F(Y)$ is the restriction of $V(f)$
on the respective space $F(X)$. Then $F$ forms covariant functor in
$Comp$, which is the subfunctor of $V$. Note that both $O$ and $OH$
are weakly normal functors.

A \emph{monad} in the category $Comp$ is a triple $\mathbb{F} =
(F,\eta,\mu)$, where $\eta:Id_{Comp}\to T$ and $\mu:T^2\to T$ are
natural transformations such that the following equalities hold: 1)
$\mu X\circ \eta F(X) = \mu X\circ F(\eta X) = \mathrm{id}_{F(X)}$;
2) $\mu X\circ \mu F(X) = \mu X\circ F(\mu X)$ \cite{EilMoore}.

The abovementioned functors generate monads. If $F$ is one of
$V,O,OH$, the identity and multiplication maps are defined as
follows. The natural transformation $\eta:Id_{{\it Comp}}\to F$ is
given by $\eta X(x)(\varphi)=\varphi(x)$ for any $x\in X$ and
$\varphi\in C(X)$, and the natural transformation $\mu: F^2\to F$
given by $\mu X(\nu)(\varphi)=\nu(\pi_{\varphi})$, where
$\pi_{\varphi}:F(X)\to \mathbb{R}$,
$\pi_{\varphi}(\lambda)=\lambda(\varphi)$.

Results on categorical and topological properties of functors $O$
and $OH$ can be found in  \cite{DOH}, \cite{OSHGEOM} \cite{OProp},
\cite{OandAR}, \cite{Conv}. In particular, we'll use the following
statements from about functors $O$ and $OH$:

\vspace*{\parindent}

{\bf Theorem 1.1.} \cite{OSHGEOM}, \cite{OProp} \emph{Functors $O$
and $OH$ are open.}

{\bf Theorem 1.2.} \cite{OSHGEOM}, \cite{OProp}, \cite{Conv}
\emph{Let $F\in\{O,OH\}$. Then 1) $\mu_F X$ is soft iff $X$ is
openly generated; 2) $F(X)\in AR$ iff $X$ is openly generated.}

{\bf Theorem 1.3.} \cite{OandAR} \emph{Let $f:X\to Y$ be open.
$O(f)$ has a degenerate fiber if and only if $f$ has}.

{\bf Theorem 1.4.} \cite{OandAR} \emph{An openly generated compactum
$X$ is $\chi$-homogeneous if and only if $O(X)$ is.}

\vspace*{\parindent}

We also note that the analogous to theorems 1.3 and 1.4 statements
hold in case of functor $OH$ and their proofs are just the same as
in case of $O$.

Finally, let us recall the definition of an $I^\tau$-fibration and
the criteria of an $I^\tau$-fibration.

A map $f:X\to Y$ is called an \emph{$I^\tau$-fibration} if it is
homeomorphic to the projection map $p_Y:I^\tau\times Y\to Y$. Note
that a map with all fibers homeomorphic to $I^\tau$ is not
necessarily an $I^\tau$-fibration (see \cite{Chig}, \cite{TorWest}
for counterexamples).

The following theorem is the well-known Torunczyk-West criterion of
a $Q$-fibration (by $Q$ we denote the Hilbert cube $[0,1]^\omega$):

\vspace*{\parindent}

{\bf Theorem 1.5.} (\cite{TorWest}) A soft mapping $f:X\to Y$ of
metric $AR$-compacta is homeomorphic to $Q$-fibration if and only if
it satisfies the condition of disjoint approximation: for any
$\varepsilon>0$ there are mappings $g_1,g_2:X\to X$ such that
$g_1(X)\cap g_2(X)=\O$, $d(g_i,id_X)<\varepsilon$, $f\circ g_i = f$.

\vspace*{\parindent}

In case of an arbitrary $\tau$, the criterion of an
$I^\tau$-fibration contains the generalization of the condition of
the Torunczyk-West theorem.

Let us give the necessary definitions first. Through $cov_\lambda
(X)$ denote the family of all coverings of cardinality $\le \lambda$
which consist of sets which are the intersections of no more than
$\lambda$ sets which are in turn are the unions of no more than
$\lambda$ co-zero sets in $X$.

A mapping $f:X\to Y$ satisfies the condition of \emph{disjoint
$\lambda$-approximation} \cite{Chig} if for any cover $\Omega\in
cov_\lambda(X)$ there exist mappings $g_1,g_2:X\to X$ with disjoint
images, $\Omega$-close to $id_X$ and with $f\circ g_i = f$.

\vspace*{\parindent}

{\bf Theorem 1.6.} (\cite{Chig}). \emph{A soft mapping $f:X\to Y$
between $AR$-compacta with fibers of weight $\le \tau$ is a trivial
$I^\tau$-fibration if and only if $f$ satisfies the condition of
disjoint $\lambda$-approximation for any $\lambda<\tau$}.

\vspace*{\parindent}

The following statement provides the sufficient condition under
which a map satisfies the condition of disjoint
$\lambda$-approximation.

\vspace*{\parindent}

{\bf Lemma 1.1} (\cite{MonInComp}). \emph{Let $f:X \to Y$ be the
limit projection $p_1$ of a $\lambda$-spectrum
$\{X_\alpha,p_\alpha,{\cal A}\}$ such that the index set ${\cal A}$
has the least element 1, all limit projections allow two disjoint
sections. Then $f$ satisfies the condition of disjoint
$\lambda$-approximation}.

\vspace*{\parindent}

\vspace*{\parindent}

{\bf 2. The main results.} For the sake of convenience, we consider
the cases $\tau = \omega$ and $\tau>\omega$ separately. Let's
consider the case of the countable $\tau$ first.

Define the metric on the space $O(X)$ for any metrizable compactum
$X$ the following way. In case $X$ is metrizable, the space $C(X)$
of all continuous functions over $X$ is separable. Choose any dense
in $C(X)$ countable set $\{\varphi_i\}_{i\in \mathbb{N}}$. We can
assume that the function $0_X$ is not in $\{\varphi_i\}_{i\in
\mathbb{N}}$. Put
$d_O(\lambda,\nu)=\sum_{i=1}^{\infty}\frac{|\lambda(\varphi_i)-\nu(\varphi_i)|}{\|\varphi_i\|\cdot
2^i}$. Then $d_O$ is an admissible metric on $O(X)$. Indeed, take
any $B_\varepsilon(\nu)=\{\lambda\in O(X)\vert
d_O(\lambda,\nu)<\varepsilon\}$. Choose number $n_0\in\mathbb{N}$
such that the inequality
$\sum_{i=n_0}^{\infty}\frac{1}{2^{i-1}}<\frac{\varepsilon}{2}$
holds. Then
$O(\nu;\varphi_1,...,\varphi_{n_0};\frac{\varepsilon}{2}\cdot
(\frac{1}{\sum_{i=n_0}^\infty \frac{1}{\|\varphi_i\|\cdot
2^i}}))\subset B_{\varepsilon}$. Hence, $d_O$ generates the topology
on $O(X)$.

Before coming to the proof of the theorem let us recall how to
extend an order-preserving functional over a single function (see
lemma 2 in \cite{OProp}).

Suppose that the set $A\subset C(X)$ is such that $0_X\in A$,
$\varphi+c_X\in A$ for any $\varphi\in A$ and $c\in\mathbb{R}$.
Consider any order-preserving functional $\nu$ on $A$ and some
function $\psi\in C(X)\backslash A$. If we want to extend $\nu$ over
the space $A\cup \{\psi+c_X\vert c\in \mathbb{R} \}$, the possible
values $\nu(\psi)$ are in the segment $[ \sup \{\nu(\varphi)\vert
\varphi\in A, \varphi\le \psi \}, \ \inf \{\nu(\varphi)\vert
\varphi\in A, \varphi\ge \psi \} ]$ and only they.

\vspace*{\parindent}

{\bf Theorem 2.1}. \emph{The mapping $\mu_{O} X$ is a $Q$-fibration
for any metrizable space $X$ which contains more than one point}.

\emph{Proof.} Assume $X$ is metrizable and not one-point. To prove
our theorem, we'll use the Torunczyk-West criterion (theorem 1.5).
It means that for any $\varepsilon>0$ we have to find two mappings
$g_1,g_2:O^2(X)\to O^2(X)$ which are both $\varepsilon$-close to
$id_{O^2(X)}$ and preserve the fibers of $\mu_O X$.

Choose some dense in $C(O(X))$ countable set $\{\Phi_i\}_{i\in
\mathbb{N}}$ that does not contain constant functions. Fix any
$\varepsilon>0$. There exists some $n_0\in \mathbb{N}$ such that
$\sum_{i=n_0}^\infty
\frac{|\Lambda(\Phi_i)-M(\Phi_i)|}{\|\Phi_i\|\cdot
2^i}<\frac{\varepsilon}{2}$ for any $\Lambda,M\in O^2(X)$.

Now for any function $\Phi_i, \ i=\overline{1,n_0}$ pick two points
$s_i,v_i\in O(X)$ such that $\Phi_i(s_i)=\sup\{\Phi_i(x) \vert x\in
O(X)\}$ and $\Phi_i(v_i)=\inf\{\Phi_i(x) \vert x\in O(X)\}$. Denote
$S_0 = \{s_i \vert i=\overline{1,n_0}\}, \ I_0 = \{v_i \vert
i=\overline{1,n_0}\}$. Since all $\Phi_i$ are continuous, for
$\varepsilon_1=\frac{\varepsilon}{b}$, where $b =
8\sum_{i=1}^{n_0}\frac{1}{\|\Phi_i\|\cdot 2^i}$, we can find
$\delta>0$ such that $B_\delta(s)\cap B_\delta(s^{'})=\O$ for any
distinct $s,s^{'}\in I_0\cup S_0$ and for any $s\in S_0\cup I_0$ and
any $x\in B_\delta(s)$ $\Phi_i(x)\in
(\Phi_i(s)-\varepsilon_1,\Phi_i(s)+\varepsilon_1)$ for all $i\in
\{1,...,n_0\}$.

Put $F_0 = \{\Phi_1,...,\Phi_{n_0}\}$. Take any $s\in I_0\cap S_0$.
Let $\Phi_{i_1},...,\Phi_{i_k}$ and $\Phi_{j_1},...,\Phi_{j_l}$ be
such that $\Phi_{i_m}(s)=\sup\{\Phi_{i_m}(x) \vert x\in O(X) \},\
m=\overline{1,k}$ and $\Phi_{j_n}(s)=\inf\{\Phi_{j_n}(x) \vert x\in
O(X)\}, n=\overline{1,l}$, where
$\{i_1,...,i_k,j_1,...,j_l\}\subset\{1,...,n_0\}$,
$\{i_1,...,i_k\}\cap\{j_1,...,j_l\}=\O$. It is clear that we can
choose distinct points $s_1,s_2\in B_\delta(s)\backslash\{s\}$ and
continuous functions
$\tilde{\Phi}_{i_1},...,\tilde{\Phi}_{i_k},\tilde{\Phi}_{j_1},...,\tilde{\Phi}_{j_l}$
such that the following conditions hold: 1)
$\tilde{\Phi}_i(x)=\Phi_i(x), \ x\in O(X)\backslash B_\delta(s), \
i\in \{i_1,...,i_k,j_1,...,j_l\}$; 2)
$\tilde{\Phi}_i(s_1)=\sup\{\tilde{\Phi}_i(x)\vert x\in O(X)\}$,
$i\in\{i_1,...,i_k\}$ and
$\tilde{\Phi}_i(s_2)=\inf\{\tilde{\Phi}_i(x)\vert x\in O(X)\}$,
$i\in\{j_1,...,j_l\}$; 3) $d(\Phi_i,\tilde{\Phi}_i)\le
\varepsilon_1$, where $i\in\{i_1,...,i_k,j_1,...,j_l\}$. We define
the new sets: $F_1 = F\cup
\{\tilde{\Phi}_{i_1},...,\tilde{\Phi}_{i_k},\tilde{\Phi}_{j_1},...,\tilde{\Phi}_{j_l}\}\backslash
\{\Phi_{i_1},...,\Phi_{i_k},\Phi_{j_1},...,\Phi_{j_l}\}$, $S_1 =
S_0\cup\{s_1\}\backslash \{s\}, \ I_1 = I_0\cup\{s_2\}\backslash
\{s\}$. We can apply the same procedure to the rest of the points of
$I_0\cap S_0$ to obtain sets $I$ and $S$ of some points of minimums
and maximums of another set of functions $F =
\{\tilde{\Phi}_1,...,\tilde{\Phi}_{n_0}\}$ with $I\cap S = \O$ and
$d(\Phi_i,\tilde{\Phi}_i)\le \varepsilon_1$ for any
$i\in\{1,...,n_0\}$. Note that $I$, for example, not necessarily
contains \emph{all} points of minimum of every function
$\tilde{\Phi}_i$, it contains only one such point; the same is about
$S$. Also we can assume that $\inf O(X)\notin S, \ \sup O(X) \notin
I$ (otherwise we would apply the reasoning from above to these
points as well).

Choose function $\Phi_0:O(X)\to \mathbb{R}$ such that
$\Phi_0(I\cup\inf O(X)) \subset \{1\}$ and $\Phi_0(S\cup\sup O(X))
\subset \{0\}$. Denote $Y = \{\pi_\varphi \ \vert \ \varphi\in
C(X)\}\cup\{\tilde{\Phi}_i+c_{O(X)} \vert \ c\in\mathbb{R},
i=\overline{1,n_0}\}$. Take any functional $M\in O^2(X)$. Let $M_0 =
\{\Lambda\in O^2(X) \ \vert \ \Lambda\vert_Y = M\vert_Y, \
\Lambda(\Phi_0) = 0\}$, and $M_1 = \{\Lambda\in O^2(X) \ \vert \
\Lambda\vert_Y = M\vert_Y, \ \Lambda(\Phi_0) = 1\}$. Due to the
choice of function $\Phi_0$, we have that $M_0\neq\O,\ M_1\neq\O$.

Let us show that the mappings $G_0,G_1:O^2(X)\to \exp O^2(X)$
defined by $G_0(M) = M_0, \ G_1(M) = M_1$ are continuous. Indeed,
take any sequence $\{M_n\}_{n\in \mathbb{N}}\subset O^2(X)$ that
converges to some $M\in O^2(X)$. We may assume that there exists $A
= \lim_{n\to\infty}(M_n)_0$. We must show that the equality $M_0 =
A$ holds. The inclusion $A\subset M_0$ is obvious. Let us show the
inclusion $M_0\subset A$ takes place. Assuming the opposite, we get
that there are $\Lambda\in M_0$ and some function $\Phi\in C(O(X))$
such that $\Lambda(\Phi)=a>\sup A(\Phi)$ or $\Lambda(\Phi)=a<\inf
A(\Phi)$ (this follows from the fact that all $(M_n)_0$ are
$O$-convex, i.e. for any $V\in O^2(X)$ with $\inf (M_n)_0\le V\le
\sup (M_n)_0$ we have $V\in (M_n)_0$, hence their limit is so (see
\cite{OandAR})). Suppose the first case holds, for instance. Note
that, since $\mu_O X$ is open, the sequence $\{\mu_O X^{-1}(\mu_O X
(M_n))\}$ converges to $\mu_O X^{-1}(\mu_O X (M))$. Hence,
$\sup\{M_n(\pi_\varphi) \ \vert \ \pi_\varphi\le\Phi, \ \varphi\in
C(X)\}$ and $\inf\{M_n(\pi_\varphi) \ \vert \ \pi_\varphi\ge\Phi, \
\varphi\in C(X)\}$ must converge to $\sup\{M(\pi_\varphi) \ \vert \
\pi_\varphi\le\Phi, \ \varphi\in C(X)\}$ and $\inf\{M(\pi_\varphi) \
\vert \ \pi_\varphi\ge\Phi, \ \varphi\in C(X)\}$ respectively.
Indeed, consider any convergent subsequence
$\{\sup\{M_{n_k}(\pi_\varphi) \ \vert \ \pi_\varphi\le\Phi, \
\varphi\in C(X)\}\}_{k\in\mathbb{N}}$ of $\{\sup\{M_n(\pi_\varphi) \
\vert \ \pi_\varphi\le\Phi, \ \varphi\in C(X)\}\}_{n\in\mathbb{N}}$,
for example (at least one such subsequence must exist!). Suppose
that its limit $s_1$ is not equal to $s = \sup\{M(\pi_\varphi) \
\vert \ \pi_\varphi\le\Phi, \ \varphi\in C(X)\}$, say $s_1>s$. Now
note that the set $\mu_O X^{-1}(\nu)$ for any $\nu\in O(X)$ consists
of all possible extensions of the functional ${\overline
\Theta}:D\to \mathbb{R}$, where $D = \{\pi_\varphi \ \vert \
\varphi\in C(X)\}$ and ${\overline
\Theta}(\pi_\varphi)=\nu(\varphi)$. Since any such extension must be
order-preserving, its possible values on $\Phi$ are in the closed
interval $[\ \sup\{{\overline \Theta}(\pi_\varphi) \ \vert \
\pi_\varphi\le\Phi, \ \varphi\in C(X)\}, \ \inf\{{\overline
\Theta}(\pi_\varphi) \ \vert \ \pi_\varphi\ge\Phi, \ \varphi\in
C(X)\} \ ]$. So, in our case we get that the possible value of any
functional from $\lim_{k\to\infty}\mu_O X^{-1}(\mu_O X(M_{n_k}))$
(again we may assume the sequence converges) cannot be less than
$s_1$ on $\Phi$, whereas functionals from $\mu_O X^{-1}(\mu_O X(M))$
are allowed to take any value up to $s$ on $\Phi$, hence $\{\mu_O
X^{-1}(\mu_O X(M_{n_k}))\}_{k\in\mathbb{N}}$ doesn't converge to
$\mu_O X^{-1}(\mu_O X(M))$, a contradiction with the openness of
$\mu_O X$. The same reasonings could be applied in the case with the
sequence $\{\inf\{M_n(\pi_\varphi) \ \vert \ \pi_\varphi\ge\Phi, \
\varphi\in C(X)\}\}_{n\in \mathbb{N}}$.

Take now any $\delta>0$. There exists $k_0\in \mathbb{N}$ such that
$\vert \sup\{M_n(\pi_\varphi) \vert \ \pi_\varphi\le\Phi, \varphi\in
C(X)\}-\sup\{M(\pi_\varphi) \vert \ \pi_\varphi\le\Phi, \varphi\in
C(X)\} \vert<\delta, \ \vert \inf\{M_n(\pi_\varphi) \vert
\pi_\varphi\ge\Phi, \varphi\in C(X)\} - \inf\{M(\pi_\varphi) \vert
\pi_\varphi\ge\Phi, \varphi\in C(X)\} \vert <\delta$ and $\vert
M(\Phi_i)-M_n(\Phi_i)\vert<\delta, \ i=\overline{1,n_0}$ for all
$n\ge k_0$. Hence, we get that $\vert \sup\{M_n(\Psi) \ \vert \
\Psi\in Y\cup \{\Phi_0\}, \ \Psi\le\Phi\}-\sup\{M(\Psi) \ \vert \
\Psi\in Y\cup \{\Phi_0\}, \ \Psi\le\Phi\}\vert <\delta$ and
$\vert\inf\{M_n(\Psi) \ \vert \ \Psi\in Y\cup \{\Phi_0\}, \
\Psi\ge\Phi\}-\inf\{M(\Psi) \ \vert \ \Psi\in Y\cup \{\Phi_0\}, \
\Psi\ge\Phi\}\vert<\delta$ for sufficiently large numbers $n$. This
means that whatever is $a=M(\Phi)\in [\sup\{M(\Psi) \ \vert \
\Psi\in Y\cup \{\Phi_0\}, \ \Psi\le\Phi\}, \ \inf\{M(\Psi) \ \vert \
\Psi\in Y\cup \{\Phi_0\}, \ \Psi\ge\Phi\}]$, we can choose $k_0\in
\mathbb{N}$ such that $[\sup\{M(\Psi) \ \vert \ \Psi\in Y\cup
\{\Phi_0\}, \ \Psi\le\Phi\}, \ \inf\{M(\Psi) \ \vert \ \Psi\in Y\cup
\{\Phi_0\}, \ \Psi\ge\Phi\}]\cap (a-\frac{a-\sup
A(\Phi)}{2},a+\frac{a-\sup A(\Phi)}{2})\neq \O$ for all $n\ge k_0$,
which means that we can obtain functionals from $(M_n)_{-}$ with
values at $\Phi$ strictly larger than $\sup A(\Phi)$, a
contradiction. Hence, the mappings $G_0,G_1$ are continuous.

Take now any $M\in O^2(X)$ and $\Lambda\in M_0$, for instance. We
have that $d_O (M,\Lambda) =
\sum_{i=1}^{\infty}\frac{|M(\Phi_i)-\Lambda(\Phi_i)|}{\|\Phi_i\|\cdot
2^i}<\sum_{i=1}^{n_0}\frac{|M(\Phi_i)-\Lambda(\Phi_i)|}{\|\Phi_i\|\cdot
2^i}+\frac{\varepsilon}{2}\le
\sum_{i=1}^{n_0}\frac{|M(\tilde{\Phi}_i)-\Lambda(\tilde{\Phi}_i)|+4\varepsilon_1}{\|\Phi_i\|\cdot
2^i}+\frac{\varepsilon}{2} =
\sum_{i=1}^{n_0}\frac{4\varepsilon_1}{\|\Phi_i\|\cdot
2^i}+\frac{\varepsilon}{2} = \varepsilon$.

Define now $g_i:O^2(X)\to O^2(X)$, $i=\overline{0,1}$  by the
formula $g_i(M) = \sup G_i$, for example. Functions $g_0, g_1$
defined that way are continuous, $\varepsilon$-close to $id_{O^2
(X)}$, with disjoint images and preserve the fibers of $\mu_O X$.
The theorem is proved.

\vspace*{\parindent}

Now we'll consider the case of $\tau>\omega$.

{\bf Theorem 2.2}. \emph{Let $w(X) = \tau>\omega$. The map $\mu_{O}
X:O^2(X)\to O(X)$ is an $I^\tau$-fibration if and only if $X$ is
openly generated and $\chi$-homogeneous.}

\emph{Proof}. {\it Sufficiency}. Suppose that $X$ is openly
generated and $\chi$-homogeneous, $w(X)=\tau>\omega$. We'll use
theorem 1.6 in combination with lemma 1.1 to prove this part of the
statement. Suppose that $\omega\le\lambda<\tau$. Represent $X$ as
the limit of a $\lambda$-system ${\cal S}=\{X_\alpha,p_\alpha,{\cal
A}\}$, where ${\cal A}$ has the minimal element $1$ and $X_1$ is a
singleton. Also we can suppose that all $p_\alpha$ are open.
Consider $Y_\alpha=O^2(X_\alpha)\times_{O(X_\alpha)}O(X)$, and by
$q_\alpha$ denote the diagonal product
$q_\alpha=(O^2(p_\alpha),\mu_{O} X)$. We obtained a $\lambda$-system
$\{Y_\alpha,q_\alpha,{\cal A}\}$ with the first limit projection
$q_1$ homeomorphic to $\mu_{O} X$. Note also that every $q_\alpha$
can be assumed soft since so is $\mu_{O} X$. We will prove that each
$q_\alpha$ allows two disjoint sections.

First let's show that the fibers of each $q_\alpha$ are infinite.
Indeed, consider any $(\Lambda,\nu)\in Y_\alpha$. Then $\mu_{O}
X_\alpha(\Lambda)=O(p_\alpha)(\nu)$. Denote $D=\{\pi_\psi \ \vert \
\psi\in C(X)\}\cup\{\Phi\circ O(p_\alpha) \ \vert \Phi\in
C(O(X_\alpha))\}$. All mappings $q_\alpha$, being soft, are
surjective. Hence, there's at least one functional
$\Theta:C(O(X))\to \mathbb{R}$ such that
$q_\alpha(\Theta)=(\Lambda,\nu)$. Then $\Theta(\pi_\psi)=\nu(\psi),
\ \Theta(\Phi\circ O(p_\alpha))=\Lambda(\Phi)$.

Our present aim is to find a function $\Phi_0\in C(O(X))$ such that
there would exist at least two distinct extensions of
$\Theta\vert_D$ on the space $D\cup \{\Phi_0+c_{O(X)} \ \vert \, \
c\in\mathbb{R}\}$.

Since $X$ is $\chi$-homogeneous and $w(X_\alpha)<w(X)$, the mapping
$p_\alpha$ doesn't have one-point fibers, and so doesn't
$O(p_\alpha)$ (theorem 1.3). Denote $S=\{\sup
O(p_\alpha)^{-1}(\lambda) \ \vert \ \lambda\in O(X_\alpha)\}$,
$I=\{\inf O(p_\alpha)^{-1}(\lambda) \ \vert \ \lambda\in
O(X_\alpha)\}$. Both these sets are closed due to the openness of
$O(p_\alpha)$ and operations $\sup,\inf:\exp O(X)\to O(X)$ being
continuous. Now define $\Phi_0\in C(O(X))$ to be a function with
$\Phi_0(S)=0$ and $\Phi_0(I)=1$. Suppose that $\Phi\circ
O(p_\alpha)\le\Phi_0$. Since $\Phi\circ O(p_\alpha)$ is constant on
the fibers of $O(p_\alpha)$, this implies $\Phi\circ O(p_\alpha)\le
0$, hence $\Phi\le 0$ and $\Theta(\Phi\circ
O(p_\alpha))=\Lambda(\Phi)\le 0$. Similarly, $\Theta(\Phi\circ
O(p_\alpha))\ge 1$ for any $\Phi\circ O(p_\alpha)\ge\Phi_0$. Now
pick any $\psi\in C(X)$ with $\pi_\psi\le \Phi_0$, for example. We
have that $\nu\in O(p_\alpha)^{-1}(\lambda)$ for some $\lambda\in
O(X_\alpha)$. Then
$\Theta(\pi_\psi)=\pi_\psi(\nu)\le\pi_\psi(\sup(O(p_\alpha)^{-1}(\lambda))\le
\Phi_0(\sup(O(p_\alpha)^{-1}(\lambda)))=0$. Similarly,
$\Theta(\pi_\psi)\ge1$ for all $\pi_\psi\ge \Phi_0$. Also, it is
obvious that $\Phi_0$ doesn't belong to $D$, hence, if we define
$\Theta(\Phi_0+c_{O(X)})=a +c$, where $a\in [0,1]$, we'll obtain an
order-preserving functional on $D\cup \{\Phi_0+c_{O(X)} \ \vert \, \
c\in\mathbb{R}\}$, which we can extend on the whole space $C(O(X))$
according to lemma 2 of \cite{OProp}. Therefore, we've shown that
$\Theta$ has at least two extensions from $D$, hence the fibers of
$q_\alpha$ are not singletons. So, for any $(\Lambda,\nu)\in
Y_\alpha$ define $g_1=\inf q_\alpha^{-1}(\Lambda,\nu)$, $g_2=\sup
q_\alpha^{-1}(\Lambda,\nu)$. The mappings $g_1,g_2$ are continuous
disjoint sections for $q_\alpha$.

Hence, the mapping $\mu_O X$ satisfies the condition of lemma 1.1,
and by theorem 1.6 it is an $I^\tau$-fibration.

\emph{Necessity}. Since $\mu_{O} X$ is an $I^\tau$-fibration, we
have that it is soft. The softness of $\mu_{O} X$ implies that $X$
is openly generated (theorem 1.2)  The space $X$ must be
$\chi$-homogeneous, since, if we suppose the opposite, we get that
$O^2(X)$ is not $\chi$-homogeneous (theorem 1.4), hence, there exist
some $\Lambda\in O^2(X)$ with $\chi(\Lambda, O^2(X))=\tau^{'}<\tau$,
and therefore $\mu_{O}X^{-1}(\mu_{O} X(\Lambda))$ is not
homeomorphic to $I^\tau$. The theorem is proved.

\vspace*{\parindent}

{\bf Note.} Proofs of theorems 2.1 and 2.2 are the same in case of
monad $\mathbb{OH}$. Note that proof of theorem 2.1 (the part of it
which concerns the choice of function $\Phi_0$) could be a bit
easier for monad $\mathbb{O}$. Indeed, the function $\Phi_0:O(X)\to
\mathbb{R}$ such that $\Phi_0(\inf O(X)) = 2\alpha$ and $\Phi_0(\sup
O(X)) = -2\alpha$, where $\alpha = \max\{\sup \Phi_i - \inf \Phi_i
\vert \ i=\overline{1,n_0}\}$ would do. In this case we can take
$M_0 = \{\Lambda\in O^2(X) \ \vert \ \Lambda\vert_Y = M\vert_Y, \
\Lambda(\Phi_0) = -\alpha\}$, and $M_1 = \{\Lambda\in O^2(X) \ \vert
\ \Lambda\vert_Y = M\vert_Y, \ \Lambda(\Phi_0) = \alpha\}$ for any
$M\in O^2(X)$, where $Y = \{\pi_\varphi \ \vert \ \varphi\in
C(X)\}\cup\{\tilde{\Phi}_i+c_{O(X)} \vert \ c\in\mathbb{R},
i=\overline{1,n_0}\}$. But in the case of $OH$ the argumentation in
proof of theorem 2.1 with $\Phi_0$, $M_0$ and $M_1$ as just
described fails.

\linespread{1}


\begin{thebibliography}{100}

\bibitem{Chig} {\it Chigogidze A.} Trivial fibrations with fibers homeomorphic to the Tychonov cube, {\it Mat. Zametki}, 1986, {\bf 39}, 747--756.

\bibitem{DOH} {\it Djabbarov G.}  Categorical properties of the functor of
weakly additive positively-homogeneous functionals. {\it Uzb. math.
jour.} {\bf 2 } 2006, p. 20-28.

\bibitem{EilMoore} {\it Eilenberg S., Moore J.} Adjoint functors and triples, {\it J.Math.}, {\bf 9}, (1965), 381--389.

\bibitem{FedChig} {\it Fedorchuk V., Chigogidze A.} Absolute retracts and infinitedimensional manifolds, - M.:Nauka, 1992, 231.

\bibitem{CovFunkt} {\it Fedorchuk V., Zarichnyii M.} Covariant functors in the category of thopological spaces, {\it Results of
Science and Thechnics. Algebra. Topology. Geometry.} {\bf 28}
Moscow.VINITI, 47--95.

\bibitem{OSHGEOM} {\it Karchevska L.} On geometric properties of
functors of positively homogeneous and semiadditive functionals. -
Preprint, 2010.

\bibitem{OProp} {\it Radul T.} On the functor of order-preserving
functionals {\it Comment.Math.Univ.Carolin}., 1998, {\bf 39},
609--615.

\bibitem{OandAR} {\it Radul T.} Topology of the space of
order-preserving functionals {\it Bull. Pol. Acad. Sci., Math.}
1999, {\bf 47}, 53--60.

\bibitem{Conv} {\it Radul T.} Convexities generated by monads. {\it
Dopov. Nats. Akad. Nauk Ukr.} Nauky 2008, 9, 27-30.

\bibitem{MonInComp} {\it Radul T., Zarichnyii M.} Monads in the category of compacta, {\it Uspekhi Mat. Nauk}, 1995, {\bf 50}, 83--106.

\bibitem{Shchepin} {\it Schepin E.} Functors and uncountable powers of compacta, {\it Uspekhi Mat. Nauk}, 1981, {\bf 36},
3--62.

\bibitem{TelZar} {\it Teleiko A., Zarichnyi M.} Categorical Topology of
Compact Hausdorff Spaces, VNTL Publishers. Lviv, 1999.

\bibitem{TorWest} {\it Torunczyk H., West J.} Fibrations and Bundles
with Hilbert Cube Manifild Fibers. - Preprint, 1982.

\bibitem{ZarMonMult} {\it Zarichnyii M.} Absolute extensors and geometry of
monad multiplication maps in the category of compacta, {\it
Mat.Sbornik}. - 1991. - {\bf 9}. - P.1261--1280 (in Russian).


\end{thebibliography}
\end{document}